\documentclass[a4paper,leqno,11pt]{amsart2000}
\usepackage{amsthm2000,amsmath2000,amscd2000,calrsfs,amsfonts,amssymb}
\usepackage{mathptmx}
\usepackage[all]{xy}
\usepackage[T1]{fontenc}

\newtheoremstyle{TH}{5pt}{5pt}{\itshape}{\parindent}{\itshape}{. ---}{ }{\thmname{#1}\thmnumber{ \bf{(#2)}}\thmnote{\normalfont{ (#3)}}}
\newtheoremstyle{DE}{5pt}{5pt}{\normalfont}{\parindent}{\itshape}{. ---}{ }{\thmname{#1}\thmnumber{ \bf{(#2)}}\thmnote{ (#3)}}
\newtheoremstyle{RE}{5pt}{5pt}{\normalfont}{\parindent}{\itshape}{. ---}{ }{\thmname{#1}\thmnumber{ \bf{(#2)}}\thmnote{ #3}\thmnote{ (#3)}}
\newtheoremstyle{PA}{}{}{\normalfont}{\parindent}{\bfseries}{ }{ }{\thmnumber{(#2)}}

\theoremstyle{TH}
\newtheorem{thm}{Theorem}[subsection]
\newtheorem{cor}[thm]{Corollary}
\newtheorem{lem}[thm]{Lemma}
\newtheorem{inlem}{Lemma}
\newtheorem{pro}[thm]{Proposition}
\newtheorem*{invertibility}{Inversion Theorem}

\theoremstyle{DE}
\newtheorem{defn}[thm]{Definition}

\theoremstyle{RE}
\newtheorem{rem}[thm]{Remark}

\newtheorem{exa}[thm]{Example}
\newtheorem{notation}[thm]{Notation}

\theoremstyle{PA}
\newtheorem{para}[thm]{Paragraph}

\newcommand{\Ref}[1]{(\ref{#1})}

\numberwithin{equation}{thm}

\newcommand{\abs}[1]{\left\vert#1\right\vert}
\newcommand{\set}[1]{\left\{#1\right\}}

\newcommand{\R}{\mathbf R}
\newcommand{\Left}{\mathbf L}
\newcommand{\C}{\mathbf C}

\newcommand{\D}{\mathbf D}

\newcommand{\Hypcoh}[1]{\mathbf{H}^{#1}}
\newcommand{\SSeqI}{{}^{I}\negmedspace E}
\newcommand{\SSeqII}{{}^{I\negthinspace I}\negmedspace E}

\newcommand{\sheaf}[1]{\mathcal{#1}}
\newcommand{\Cat}[1]{\mathfrak{#1}}

\newcommand{\Ox}[1]{\mathcal{O}_{#1}}

\newcommand{\Spec}[1]{\mathrm{Spec}(#1)}
\newcommand{\Pico}[1]{\mathrm{Pic}^{\circ}\left(#1\right)}
\newcommand{\Jac}[1]{{\mathrm{J}(#1)}}
\newcommand{\Proj}[2][]{\mathbf{P}^{\negthinspace#2}_{\negthinspace #1}} 
\newcommand{\Hom}{\mathrm{Hom}}

\newcommand{\proj}[1]{\mathrm{pr}_{#1}}
\newcommand{\maps}{\colon}

\newcommand{\id}[1]{1_{#1}}

\newcommand{\complex}[1]{#1^{\bullet}}

\DeclareMathOperator{\rank}{rk}
\DeclareMathOperator{\coker}{coker}
\newcommand{\ch}{\mathrm{ch}}

\newcommand{\Tensor}{\overset{\Left}{\otimes}}

\makeatother
\providecommand{\specseq}[9]{%
    \xymatrix{%
        #1\quad=&{}\ar@{-}'[]+<0pt,10pt>'[dd]+<0pt,15pt>_<<{q}'[ddrrr]+<-10pt,15pt>_>>{p}&%
        {#2}\ar@{#6}[r]^-{#7}&{#3}\\&&%
        {#4}\ar@{#8}[r]^-{#9}&{#5}\\&&&&%
    }%
} \makeatletter

\newcommand{\higgs}[3][X]{(\sheaf{#2}\xrightarrow{#3}\sheaf{#2}\otimes\omega_{#1})}
\newcommand{\Higgs}[2][]{\mathsf{#2}_{#1}}
\newcommand{\Fourier}[1]{\mathbf{F}_{#1}}
\newcommand{\Mukai}[1]{\mathbf{M}_{#1}}
\newcommand{\ITF}[2]{{\Phi}^{#1}_{#2}}
\newcommand{\TFT}[1]{\mathbf{TF}(#1)}
\newcommand{\THC}[1]{\complex{\sheaf{H}}(\mathsf{#1})}

\newcommand{\UTensor}{\overset{\Left}{\boxtimes}}

\begin{document}


\title[]{A Fourier transform for Higgs bundles}
\author{J. Bonsdorff}
\address{Mathematical Institute\\
    24$-$29 St Giles\\
    Oxford OX1 3LB\\
    United Kingdom}
\email{Juhani.Bonsdorff@merton.ox.ac.uk}

\thanks{The author was supported by The Finnish Cultural Foundation, The Väisälä
Foundation, and The Scatcherd European Scholarship}


\begin{abstract}
We  define a Fourier-Mukai transform for Higgs bundles on smooth
curves (over  $\C$ or another algebraically closed field) and
study its properties. The transform of a stable degree-$0$ Higgs
bundle is an algebraic vector bundle on  the cotangent bundle of
the Jacobian of the curve. We show that the transform admits a
natural extension to an algebraic vector bundle over projective
compactification of the base. The main result is that the
original Higgs bundle can be reconstructed from this extension.
\end{abstract}

\maketitle


\tableofcontents


\section{Introduction}
Higgs bundles on a compact Riemann surface $X$ are pairs
$(E,\theta)$ consisting of a holomorphic vector bundle $E$ and a
holomorphic one-form $\theta$ with values in $\mathrm{End}(E)$ on
$X$. They originated essentially simultaneously in  Nigel
Hitchin's study \cite{hitchin:self-duality} of dimensionally
reduced self-duality equations of Yang-Mills gauge theory and,
over general Kähler manifolds, in Carlos Simpson's work
\cite{simpson:constructing} on Hodge theory.

To explain Hitchin's point of view, we consider  solutions of the
self-dual Yang-Mills equations  on $\R^4$ that are invariant
under translations in one or more directions in $\R^4$.
Invariance in one direction reduces the SDYM equations to the
Bogomolnyi equations describing magnetic monopoles in $\R^3$,
while invariance in three directions leads to Nahm's equation on
$\R$. Invariance in two directions leads to the conformally
invariant Hitchin's equations (or {\em Higgs bundle equations}),
which the conformal invariance allows to be considered on compact
Riemann surfaces. A solution to Hitchin's equations has an
interpretation as a triple $(E,\theta,h)$ with $E$  a holomorphic
 vector bundle, $\theta$ a holomorphic section of
$\mathrm{End}(E)\otimes\omega_X$ and $h$ a Hermitian metric on
$E$, satisfying $F+[\theta,\theta^* ]=0$, where $F$ is the
curvature of the connection determined by the metric. It was
shown in Hitchin \cite{hitchin:self-duality} for rank-$2$ bundles
on Riemann surfaces, and in Simpson \cite{simpson:constructing}
in general, that a pair $(E,\theta)$ admits  a unique such metric
precisely when $E$ has vanishing Chern classes and the pair
$(E,\theta)$ is {\em stable} in a sense which generalises the
usual stability for vector bundles. For an excellent overview of
Simpson's viewpoint of  non-Abelian cohomology, see Simpson
\cite{simpson:hodge}.

An important class of  transforms in Yang-Mills theory, including
the ADHM construction and the Fourier transform for instantons
(Donaldson-Kronheimer \cite{Donaldson-Kronheimer}) and the Nahm
transform for monopoles, is based on using the kernel of the Dirac
operator coupled to the (twisted) connection.  A recent work of
Marcos Jardim \cite{jardim:thesis,jardim:constructing} uses a
version of the Nahm transform to establish a link between
singular Higgs bundles on a $2$-torus and doubly-periodic
instantons. Our goal is to generalise this work to Riemann
surfaces of genus $g\geq 2$. In this paper we shall consider the
purely holomorphic aspect of the transform; we plan to return to
the properly gauge-theoretic questions in a future paper.

The holomorphic side of the Nahm transform is captured by the
(generalised) Fourier-Mukai transform, which originated in the
work of Shigeru Mukai \cite{mukai:duality} on Abelian varieties.
Let $X$ be a complex torus and $\Hat{X}$ its dual, and let
$\mathbf{D}(X)$ and $\mathbf{D}(\Hat{X})$ denote the {\em derived
categories} of the categories of coherent sheaves on $X$ and
$\Hat{X}$ respectively. Using the {\em Poincar\'{e} sheaf}
$\sheaf{P}$ on $X\times\Hat{X}$, Mukai defined a functor
$M\maps\mathbf{D}(X)\to\mathbf{D}(\Hat{X})$ by
$$
M(\bullet)=\R\proj{\Hat{X}*}(\proj{X}^*(\bullet)\otimes\sheaf{P}),
$$
and showed that it is a category equivalence. This construction
can be generalised to any varieties $X$ and $Y$ together with a
sheaf $\sheaf{P}$ on $X\times Y$. The properties of these
generalisations have been studied by A. Bondal and D. Orlov
\cite{bondal-orlov:semiorthogonal,bondal-orlov:reconstructing},
A. Maciocia \cite{maciocia:generalized}, T. Bridgeland
\cite{bridgeland:thesis,bridgeland:equivalences}, and others.

We interpret the endomorphism-valued one-form $\theta$ as a bundle
map, making a Higgs bundle $(E,\theta)$ into a sheaf complex
$\sheaf{E}\to \sheaf{E}\otimes\omega_X$, where $\omega_X$ is the
canonical sheaf of $X$. Hence a Higgs bundle gives us an object
of the derived category $\mathbf{D}(X)$, and we can use the
general machinery of Fourier-Mukai transforms. The key
observation is that it is necessary to consider relative
transforms of families of Higgs bundles twisted by adding a
scalar term to the Higgs field $\theta$. Our transform produces
sheaves on $\Jac{X}\times H^0(X,\omega_X)$, where $\Jac{X}$ is the
Jacobian of $X$. This base can be identified with the cotangent
bundle of the Jacobian.

While the motivation for the present work comes from differential
and complex analytic geometry, we are actually working within the
framework of algebraic geometry, noting that on an algebraic
curve the Higgs bundle data is purely algebraic. Translation
between these categories is provided by Serre's GAGA \cite{GAGA}.
Our approach has the advantage that all  constructions are
automatically algebraic (or holomorphic), while the fact that we
are dealing with rather high-dimensional base spaces would make
some of the analytic techniques of Jardim hard to use.

The first part of this paper develops the  machinery of
generalised Fourier-Mukai transforms.  While some of the material
presented in section 2 cannot be found in the literature, it is
mostly well known. The main new contributions are the definition
of a relative Fourier transform for curves and the reduction of
it to the original Mukai transform.

The transform for Higgs bundles is developed in section 3. The
first interesting application  is that our Fourier transform
takes stable Higgs bundles of degree zero to vector bundles on
$\Jac{X}\times H^0(X,\omega_X)$. Our approach has the advantage
of giving directly an algebraic (holomorphic) extension of this
bundle over the projective completion
$\Jac{X}\times\Proj{}(H^0(X,\omega_X)\oplus \C)$ of the base
space, without a need to separately compactify a bundle on
$\Jac{X}\times H^0(X,\omega_X)$. Denote this extension of the
transform of a Higgs bundle $\Higgs{E}=(E,\theta)$ by
$\TFT{\Higgs{E}}$. The main theorem of this paper is the
following:
\begin{invertibility}[\ref{thm:invertibility}]
Let $\Higgs{E}$ and $\Higgs{F}$ be two Higgs bundles on a curve $X$ of genus $g\geq 2$.
If $\TFT{\Higgs{E}}\cong\TFT{\Higgs{F}}$, then
$\Higgs{E}\cong\Higgs{F}$ {\em as Higgs bundles}.
\end{invertibility}
We in fact prove this  theorem by exhibiting a procedure for
recovering a Higgs bundle from its transform. Furthermore, it
follows easily from the theorem that the transform functor is in
fact fully faithful.

\vspace{17pt}

\noindent {\bf Acknowledgements.} The original idea of developing
a Fourier transform for Higgs bundles is due to Nigel Hitchin
\cite{hitchin:dirac}. I am deeply grateful to him for generous
comments, support and encouragement.

\subsection*{Notation and conventions}
Unless otherwise specified, all rings and algebras are
commutative and unital. We fix an algebraically closed field $k$.
All schemes are assumed to be of finite type over  $k$. All
morphisms are $k$-morphisms and all products are products over
$\Spec{k}$ unless stated otherwise. A curve always means  a
smooth irreducible complete (i.e., projective) curve over $k$. If
$\sheaf{F}$ is an $\Ox{X}$-module, $\sheaf{F}^{\vee}$ denotes its
dual $\sheaf{H}om_{\Ox{X}}(\sheaf{F},\Ox{X})$. The canonical
sheaf of a curve $X$ is denoted by $\omega_X$.

$\D(X)$ denotes the derived category of the (Abelian) category of
$\Ox{X}$-modules. $\D^-(X)$ (resp. $\D^+(X)$, resp. $\D^b(X)$) is
the full subcategory of objects cohomologically bounded above
(resp. bounded below, resp. bounded). $\D_{qcoh}(X)$ and
$\D_{coh}(X)$ are the full subcategories of objects with
quasi-coherent and coherent cohomology objects, respectively.
These superscripts and subscripts can be combined in the obvious
way. The category of $\Ox{X}$-modules is denoted by
$\Cat{Mod}(X)$, and $\Cat{QCoh}(X)$ is the thick subcategory of
quasi-coherent sheaves.

A commutative square
$$
\begin{CD}
Z @>v>> X\\
@VgVV   @VVfV\\
Y @>>u> S
\end{CD}
$$
is called {\em Cartesian} if the mapping $(v,g)_S\maps Z\to
X\times_S Y$ is an isomorphism. We denote {\em canonical}
isomorphisms often by "$=$".

\section{Fourier-Mukai transforms}
We develop here the general machinery of Fourier-Mukai transforms
that will be necessary for our application to Higgs bundles.

We will be using the theory of derived categories;  the main
reference to derived categories in algebraic geometry remains
Hartshorne's seminar \cite{hartshorne:residues} on
Grothen{-}dieck's duality theory. Further references include
Gelfand-Manin \cite{gelfand-manin:methods}, Kashiwara-Shapira
\cite{kashiwara-shapira:sheaves} and Weibel
\cite{weibel:introduction}. For a nice informal introduction, see
Illusie \cite{illusie:categories} or the introduction of Verdier's
thesis \cite{verdier:thesis}.

\subsection{A base change result}
\begin{para}
Consider the following diagram of schemes (here not necessarily
of finite type over a field):
\begin{equation*}
    \xymatrix{
        & X_1\times_S X_2   \ar[ld]_{p_1}
                            \ar[dd]_{f}
                            \ar[rd]^{p_2}               \\
        X_1\ar[dd]_{f_1}    &&          X_2\ar[dd]^{f_2}\\
        & Y_1\times_S Y_2   \ar[ld]_{q_1}
                            \ar[rd]^{q_2}               \\
        Y_1\ar[rd]          &&          Y_2\ar[ld]      \\
        & S
    }
\end{equation*}
with $f=f_1\times_S f_2$. Recall  the external tensor product over
$S$ of an $\Ox{X_1}$-module  $\sheaf{F}_1$ and an
$\Ox{X_2}$-module $\sheaf{F}_2$:
$$
\sheaf{F}_1\boxtimes_S \sheaf{F}_2 =
\left({p_1}^*\sheaf{F}_1\right)\,\otimes_{\Ox{X_1\times_S X_2}}\,
\left({p_2}^*\sheaf{F}_2\right).
$$
We get the corresponding left-derived bifunctor
$$
(\bullet)\UTensor_S(\bullet)\maps\D^{-}(X_1)\times\D^{-}(X_2)\to
\D^{-}(X_1\times_S X_2).
$$
\end{para}

The following theorem should be part of folklore; we include a
proof of it for the lack of a suitable reference. It is
essentially the derived-category version of a part of
Grothendieck's theory of "global hypertor functors" (EGA III
\cite{EGAIII}, \S 6).

\begin{thm}[Künneth formula]\label{thm:kunneth}
For $i=1,2$ let $\sheaf{F}_i$ be an object of $\D_{qcoh}^-(X_i)$.
Assume that the schemes are Noetherian and of finite dimension,
and that the $f_i$ are separated. Then
$$\left(\R f_{1*}\sheaf{F}_1 \right)\UTensor_S \left(\R f_{2*}\sheaf{F}_2\right) =
\R f_*\left(\sheaf{F}_1\UTensor_S\sheaf{F}_2\right)$$ if either
$\sheaf{F}_1$ or $\sheaf{F}_2$ is quasi-isomorphic to a  complex
of $S$-flat sheaves. This is true in particular  if either $X_1$
or $X_2$ is flat over $S$.
\end{thm}

\begin{proof}
The Noetherian and dimensional hypotheses guarantee that the
derived direct images are defined for complexes not bounded
below. There are natural "adjunction" maps $\id{}\to\R f_*\Left
f^*$ giving
$$
(\R f_{1*}\sheaf{F}_1) \UTensor_S (\R f_{2*}\sheaf{F}_2) \to \R
f_* \Left f^*\left((\R f_{1*}\sheaf{F}_1) \UTensor_S (\R
f_{2*}\sheaf{F}_2)\right).
$$
Notice that
\begin{align*}
\Left f^*\left((\R f_{1*}\sheaf{F}_1)\UTensor_S (\R
f_{2*}\sheaf{F}_2)\right) &=
(\Left f^*\Left q_1^*\R f_{1*}\sheaf{F}_1)\Tensor (\Left f^*\Left q_2^*\R f_{2*}\sheaf{F}_2)\\
&=(\Left p_1^*\Left f_1^*\R f_{1*}\sheaf{F}_1) \Tensor (\Left
p_2^*\Left f_2^*\R f_{2*}\sheaf{F}_2).
\end{align*}
Now the adjunctions $\Left f_i^*\R f_{i*}\to\id{}$ give a natural
map
\begin{align*}
(\Left p_1^*\Left f_1^*\R f_{1*}\sheaf{F}_1) \Tensor (\Left
p_2^*\Left f_2^*\R f_{2*}\sheaf{F}_2) &\to (\Left
p_1^*\sheaf{F}_1)\Tensor (\Left p_2^*\sheaf{F}_2)\\
&=\sheaf{F}_1\UTensor_S\sheaf{F}_2.
\end{align*}
Composing  gives us a natural transformation
$$
(\R f_{1*}\sheaf{F}_1) \UTensor_S (\R f_{2*}\sheaf{F}_2) \to \R
f_*\left(\sheaf{F}_1\UTensor_S\sheaf{F}_2\right).
$$
Whether this is an isomorphism is a local question; hence we may
assume that $S=\Spec{A}$ and $Y_i=\Spec{B_i}$ for $i=1,2$.

Suppose $\sheaf{F}_1$ is quasi-isomorphic to a complex of
$S$-flat sheaves; replace $\sheaf{F}_1$ with this flat
resolution. Then $\sheaf{F}_1\UTensor_S\sheaf{F}_2=
\sheaf{F}_1\boxtimes_S\sheaf{F}_2$.

For $i=1,2$ let $\mathfrak{U}_i=(U_{i,\alpha})_{\alpha}$ be a
finite affine open cover of $X_i$. Let $\mathfrak{U}$ denote the
open affine cover $(U_{1,\alpha}\times_S
U_{2,\beta})_{\alpha,\beta}$ of $X_1\times_S X_2$. Notice that in
all these covers arbitrary intersections  of the covering sets
are affine. Let $\complex{\check{C}}(\mathfrak{U}_i,\sheaf{F}_i)$
denote the simple complex associated to the \v{C}ech double
complex of $\sheaf{F}_i$ with respect to $\mathfrak{U}_i$.
Similarly, let
$\complex{\check{C}}(\mathfrak{U},\sheaf{F}_1\boxtimes_S\sheaf{F}_2)$
be the \v{C}ech complex with respect to $\mathfrak{U}$.

Now $\R\Gamma(X_i,\sheaf{F}_i)$ is quasi-isomorphic to
$\complex{\check{C}}(\mathfrak{U}_i,\sheaf{F}_i)$, and hence $\R
f_{i*}\sheaf{F}_i$ is quasi-isomorphic to
$\complex{\check{C}}(\mathfrak{U}_i,\sheaf{F}_i)^{\sim}$. But the
sheaves of these complexes are $S$-flat by construction, whence
$$
(\R f_{1*}\sheaf{F}_1)\UTensor_S (\R f_{2*}\sheaf{F}_2) =
\left(\complex{\check{C}}(\mathfrak{U}_1,\sheaf{F}_1)\otimes_A
\complex{\check{C}}(\mathfrak{U}_2,\sheaf{F}_2)\right)^{\sim}.
$$
Similarly
$$
\R f_*\left(\sheaf{F}_1\UTensor_S\sheaf{F}_2\right) =
\left(\complex{\check{C}}(\mathfrak{U},\sheaf{F}_1\boxtimes_S\sheaf{F}_2)\right)^{\sim}.
$$
Hence we are reduced to showing that the complex
$\complex{\check{C}}(\mathfrak{U}_1,\sheaf{F}_1)\otimes_A
\complex{\check{C}}(\mathfrak{U}_2,\sheaf{F}_2)$ is
quasi-isomorphic to
$\complex{\check{C}}(\mathfrak{U},\sheaf{F}_1\boxtimes_S\sheaf{F}_2)$.
But  this is showed in the proof of (6.7.6) of EGA III
\cite{EGAIII}.
\end{proof}

\begin{rem}
If one wants to avoid the Noetherian hypothesis in the theorem,
one could work with objects $\sheaf{F}_i$ of
$\D^-(\Cat{QCoh}(X_i))$ and require the $f_i$ to be quasi-compact.
This is essentially the viewpoint of EGA III.
\end{rem}

\begin{cor}\label{cor:basechange}
Let $f\maps X\to S$ and $g\maps T\to S$ be  morphisms of
finite-dimensional Noetherian schemes. Let $f'\maps X\times_ST\to
T$ and $g'\maps X\times_ST\to X$ be the projections, and let
$\sheaf{F}$ belong to $\D^-_{qcoh}(X)$.
\begin{enumerate}
\item If $\sheaf{F}$ is quasi-isomorphic to a complex of $S$-flat
sheaves (in particular, if $f$ is flat), then
$$
\Left g^* \R f_* \sheaf{F} = \R f'_{*} \,\Left {g'}^* \sheaf{F}.
$$
\item If $g$ is flat, then
$$
g^*\R f_* \sheaf{F} = \R f'_{*} {g'}^*\sheaf{F}.
$$
\end{enumerate}
\end{cor}
\begin{proof}
Apply the Künneth formula with $X_1=X$, $Y_1=S$, $f_1=f$,
$X_2=Y_2=T$, $f_2=\id{T}$, $\sheaf{F}_1=\sheaf{F}$ and
$\sheaf{F}_2=\Ox{T}$.
\end{proof}

\subsection{Integral transforms}
\begin{defn}
Let $S$ be a separated $k$-scheme and let $X$ and $Y$ be flat
$S$-schemes. If $\sheaf{P}$ is an object of $\D_{coh}^b(X\times_S
Y)$, the {\em relative integral transform} defined by $\sheaf{P}$
is the functor $\ITF{\sheaf{P}}{X\to Y/S}\maps\D^+(X)\to\D^+(Y)$
given by
$$
\ITF{\sheaf{P}}{X\to
Y/S}(\bullet)=\R{\proj{2}}_*(\proj{1}^*(\bullet)\Tensor\sheaf{P}),
$$
where $\proj{1}$ and $\proj{2}$ are the canonical projections of
$X\times_S Y$. When $S=\Spec{k}$ we call the transform the {\em
absolute integral transform} and denote it by
$\ITF{\sheaf{P}}{X\to Y}$.
\end{defn}

\begin{pro}\label{pro:rel->abs}
Let $i\maps X\times_S Y\to X\times_k Y$ be the morphism
$(\proj{1},\proj{2})_k$. Then $\R i_*=i_*$ and
$$
\ITF{\sheaf{P}}{X\to Y/S}(\bullet) = \ITF{i_*\sheaf{P}}{X\to Y}(\bullet).
$$
\end{pro}
\begin{proof}
We have the commutative diagram
\begin{equation*}
\xymatrix{
{}  &   X\times_S Y \ar[dl]_{\proj{1}}
                    \ar[dr]^{\proj{2}}
                    \ar[dd]_i            \\
X   &   { }                  &       Y   \\
{ } &   X\times Y.  \ar[ul]^p
                    \ar[ur]_q
}
\end{equation*}
Notice that because both $\proj{1}$ and $p$ are flat morphisms, we
have
$$
    \proj{1}^*=\Left\proj{1}^*=\Left(i^*\circ p^*)=\Left i^*\circ\Left p^* =
    \Left i^*\circ p^*.
$$
Using this and the projection formula, we have
\begin{align*}
\ITF{\sheaf{P}}{X\to Y/S}(\bullet)
        &= \R{\proj{2}}_*(\proj{1}^*(\bullet)\Tensor\sheaf{P})\\
        &= \R q_*\R i_* (\Left  i^*
           (p^*(\bullet))\Tensor\sheaf{P})\\
        &= \R q_* (p^*(\bullet)\Tensor\R i_*\sheaf{P}).
\end{align*}
But $i$ fits in a Cartesian square
\begin{equation*}
\begin{CD}
X\times_S Y @>{i}>> X\times_k   Y            \\
@VVV                @VVV                     \\
S           @>>{\Delta_{S/k}}>  S\times_k S. \\
\end{CD}
\end{equation*}
As $S/k$ is separated, $\Delta_{S/k}$ is a closed immersion, and
consequently so is $i$. In particular, $i_*$ is an exact functor
and therefore equal to $\R i_*$. Hence
\begin{align*}
\ITF{\sheaf{P}}{X\to Y/S}(\bullet)
                 = \R q_* (p^*(\bullet)\Tensor i_*\sheaf{P})
                 = \ITF{i_*\sheaf{P}}{X\to Y}(\bullet)
\end{align*}
as claimed.
\end{proof}
\begin{rem}
We cannot avoid using the derived tensor product in the above
result, even if $\sheaf{P}$ is a locally free sheaf, because
$i_*\sheaf{P}$ is not flat in general. However, as $i$ is proper,
$i_*\sheaf{P}$ belongs always to $\D^b_{coh}(X\times Y)$.
\end{rem}

\begin{para}\label{para:Fs(k(x))}
For flat $S$-schemes $X$ and $Y$ and for $x\in X$, let $Y_x$
denote the fibre $\proj{1}^{-1}(x)$, where $\proj{1}\maps
X\times_S Y\to X$ is the canonical projection. We have then a
commutative diagram
$$
\begin{CD}
Y_x         @>j>>    X\times_S Y     @>{\proj{2}}>>    Y      \\
@VVV                @V{\proj{1}}VV                    @VVV   \\
\kappa(x)   @>>>    X               @>>>    S,     \\
\end{CD}
$$
in which all squares are Cartesian. Let $i$ denote the composition
of the top arrows. For an object $\sheaf{F}$ of
$\D_{coh}^b(X\times_S Y)$ (resp.  $\D_{coh}^b(Y)$), we denote by
$\sheaf{F}_x$  the "restriction" $\Left j^*\sheaf{F}$ (resp.
$\Left i^*\sheaf{F}$) to $Y_x$. For complexes of locally free
sheaves these are just ordinary restrictions to $Y_x$. If
$\sheaf{P}$ is a locally free sheaf on $X\times_S Y$, then for
each $x\in X$
$$
    \ITF{\sheaf{P}}{X\to Y/S}(k(x)) = i_*\sheaf{P}_x,
$$
where $k(x)$ is the skyscraper sheaf $k$ at $x$. Indeed, consider
the commutative diagram above: the claim follows from flat base
change around the left-hand square and the projection formula
applied to $j$. Notice that $i_*$ is exact.
\end{para}
\begin{exa}\label{exa:relative_mukai}
Let $X$ be an Abelian variety, $\Hat{X}$ its dual, and let $S$ be
a separated scheme. Let $\sheaf{P}$ be the {\em Poincaré sheaf} on
$X\times\Hat{X}$, normalised so that both
$\sheaf{P}|_{X\times\set{0}}$ and
$\sheaf{P}_{\set{0}\times\Hat{X}}$ are the trivial line bundles.
Denote by $\sheaf{P}_S$ the pull-back of this Poincar\'{e} sheaf
to $X\times \Hat{X}\times S = (X\times S)\times_S(\Hat{X}\times
S)$. The {\em relative Mukai transform functor}
$\Mukai{S}\maps\D^b_{coh}(X\times S)\to\D^b_{coh}(\Hat{X}\times
S)$ is the relative integral transform functor
$\ITF{\sheaf{P}_S}{(X\times S)\to (\Hat{X}\times S)/S}$. If
$S=\Spec{k}$, we denote the transform by $\Mukai{}$.
\end{exa}

The following theorem of Mukai plays a crucial role in the proof
of our invertibility result \Ref{thm:invertibility}.

\begin{thm}\label{thm:relative_mukai_equivalence}
If $S$ is a smooth projective variety, then the relative Mukai
transform $\Mukai{S}$ is an equivalence of  categories from
$\D^b_{coh}(X\times S)$ to  $\D^b_{coh}(\Hat{X}\times S)$.
\end{thm}
\begin{proof}
See Mukai \cite{mukai:fourier}. The proof is a generalisation of
Mukai's original proof of this result for the absolute transform
$\Mukai{}$ in \cite{mukai:duality}.
\end{proof}

\begin{pro}\label{pro:ITF-basechange}
Let $X$ and $Y$ be flat $S$-schemes and $\sheaf{P}$ an object of
$\D^b_{coh}(X\times_S Y)$. Let $u\maps T\to S$ be a morphism of
schemes. Let $i_X\maps X_{(T)}\to X$, $i_Y\maps Y_{(T)}\to Y$,
and $j\maps (X\times_S Y)_{(T)}=X_{(T)}\times_T Y_{(T)}\to
X\times_S Y$ be the canonical projections. Then
$$
\Left i_Y^* \circ \ITF{\sheaf{P}}{X\to Y/S}
=
\ITF{\Left j^*\sheaf{P}}{X_{(T)}\to Y_{(T)}/T}\circ\Left i_X^*.
$$
Moreover, if $u$ is a flat morphism, then all derived pull-backs
above can be replaced with normal pull-backs.
\end{pro}
\begin{proof}
Consider the commutative diagram
\begin{equation*}
\xymatrix@R=20pt@!C=30pt{
X_{(T)}\times_T Y_{(T)}
    \ar[ddd]_q
    \ar[rrr]^p
    \ar[rd]_j
&&& X_{(T)}
    \ar[dl]^{i_X}
    \ar[ddd]                    \\
& X\times_S Y
    \ar[d]_{\proj{2}}
    \ar[r]^{\proj{1}}
&X  \ar[d]                      \\
& Y \ar[r]
& S                             \\
Y_{(T)}
    \ar[ur]^{i_Y}
    \ar[rrr]
&&& T
    \ar[ul]_u                   \\
}
\end{equation*}
It is immediate that all squares are Cartesian. If $u$ is flat,
then so are $i_X$, $i_Y$ and $j$; this proves the claim about
replacing derived pull-backs. Since in any case $X/S$ is flat,
$\proj{2}$ is also flat. So by \Ref{cor:basechange}  we can do a
base change around the leftmost square. We get
\begin{align*}
\Left i_Y^* \ITF{\sheaf{P}}{X\to Y/S}(\bullet)
    &= \Left i_Y^*\R\proj{2*}\left(\proj{1}^*(\bullet)\Tensor\sheaf{P}\right)
    \\
    &= \R q_* \Left j^*\left(\proj{1}^*(\bullet)\Tensor\sheaf{P}\right)
    \\
    &= \R q_* \left(\Left j^*(\proj{1}^*(\bullet))\Tensor\Left j^*\sheaf{P}\right)
    \\
    &= \R q_*\left( p^*\Left i_{X}^*(\bullet)\Tensor\Left j^*\sheaf{P}\right)
    =\ITF{\Left j^*\sheaf{P}}{X_{(T)}\to Y_{(T)}/T}(\Left
    i_X^*(\bullet)).
\end{align*}
\end{proof}

\begin{pro}\label{pro:RGamma}
Let $X$ and $Y$ be flat $S$-schemes and $\sheaf{P}$ an object of
 $\D^b_{coh}(X\times_S Y)$. Then
$$
\R\Gamma(Y,\ITF{\sheaf{P}}{X\to Y/S}(\sheaf{E}))=
\R\Gamma(X,\sheaf{E}\Tensor\R\proj{1*}\sheaf{P}).
$$
\end{pro}
\begin{proof}
We simply use the composition property of derived functors  and
the projection formula:
\begin{align*}
\R\Gamma(Y,\ITF{\sheaf{P}}{X\to Y/S}(\sheaf{E})) &=
\R\Gamma(Y,\R\proj{2*}(\proj{1}^*\sheaf{E}\Tensor\sheaf{P}))
        &\text{(by definition)}\\
&= \R\Gamma(X\times_SY,\proj{1}^*\sheaf{E}\Tensor\sheaf{P})
        &\text{(composition)}\\
&= \R\Gamma(X,\R\proj{1*}(\proj{1}^*\sheaf{E}\Tensor\sheaf{P}))
        &\text{(composition)}\\
&=  \R\Gamma(X,\sheaf{E}\Tensor\R\proj{1*}\sheaf{P})
        &\text{(by projection formula)}.
\end{align*}
\end{proof}

\subsection{WIT complexes}
\begin{notation}\label{notation:WIT}
Let $X$ and $Y$ be proper flat  $S$ schemes. We fix a locally free
sheaf $\sheaf{P}$ on $X\times_S Y$, and denote by $F_S$ the
relative integral transform functor $\ITF{\sheaf{P}}{X\to
Y/S}\maps\D^b_{coh}(X)\to\D^b_{coh}(Y)$. We leave it to the
reader to generalise the results of this subsection to a more
general setting.
\end{notation}

\begin{defn}
We say that an object $\sheaf{E}$ of $\D_{coh}^b(X)$ is a {\em
$WIT_{\sheaf{P}}(n)$-complex}\footnote{Following Mukai, "WIT"
stands for "weak index theorem".} if $H^p(F_S(\sheaf{E}))=0$ for
all $p\neq n$. If  $\sheaf{P}$ is clear from the context, we
shall omit the explicit reference to it. An object of
$\D_{coh}^b(X)$ is a WIT-{\em complex} if it is a $WIT(n)$-complex
for some $n$.

If $\sheaf{E}$ is a $WIT(n)$-complex on $X$, the (coherent) sheaf
$H^n(F_S(\sheaf{E}))$ on $Y$ is called the {\em integral
transform of $\sheaf{E}$}, and is denoted by
$\widehat{\sheaf{E}}$.
\end{defn}

\begin{defn}
We say that an object $\sheaf{E}$ of $\D_{coh}^b(X)$ is an
$IT_{\sheaf{P}}(n)$-complex\footnote{ "IT" stands for "Index
theorem".} if for each (closed) point $y\in Y$ and each $p\neq n$
we have
$$
\Hypcoh{p}(X_y,\sheaf{E}_y\otimes\sheaf{P}_y)=0,
$$
where we are using the notation of \Ref{para:Fs(k(x))} for
$\sheaf{E}_y$, $\sheaf{P}_y$ and $X_y$.
\end{defn}

\begin{lem}\label{lem:basechange}
Let $f\maps X\to Y$ be a proper morphism of (Noetherian) schemes
and let $\mathcal{E}$ be an object of $\D_{coh}^b(X)$ which has a
$Y$-flat resolution. Let $y\in Y$. Then:
\begin{enumerate}
  \item if the natural map $\varphi^p(y)\maps\mathbf{R}^pf_*(\mathcal{E}) \otimes \kappa(y)\to
  \mathbf{H}^p(X_y,\mathcal{E}_y)$ is surjective, then it is an
  isomorphism.
  \item If $\varphi^p(y)$ is an isomorphism, then $\varphi^{p-1}$ is also an
  isomorphism if and only if
  $\mathbf{R}^pf_*(\mathcal{E})$is free in a open
  neighbourhood of $y$.
\end{enumerate}
\end{lem}

\begin{proof}
This follows from  EGA III \cite{EGAIII} \S7. However,  that part
of EGA can be somewhat hard to read; one could also follow the
simpler proof of Hartshorne \cite{AG} Theorem III.12.11, making
the fairly minor and obvious adjustments for hypercohomology.
\end{proof}

\begin{pro}\label{pro:IT}
Let $\sheaf{E}$ be an $IT(n)$ complex. Then $\sheaf{E}$ is a
$WIT(n)$-complex, and $\widehat{\sheaf{E}}$ is locally free on
$Y$.
\end{pro}

\begin{proof}
Our schemes are Jacobson, and so it suffices to restrict our
attention to closed points. Since $\proj{2}$ is flat,
$\proj{1}^*\sheaf{E}$ is quasi-isomorphic to a complex of sheaves
flat over $Y$. Moreover, $X$ is proper over $S$, and so
$\proj{2}$ is a proper morphism. We are then in position to use
\Ref{lem:basechange}. Let $y\in Y$ be a closed point. Now
$$
(\proj{1}^*\mathcal{E}\otimes\sheaf{P})_y \cong
\mathcal{E}_y\otimes\mathcal{P}_y
$$
on $(X\times_S  Y)_{y}= X_y$. Hence by hypothesis the natural map
$$
\varphi^p(y)\maps\mathbf{R}^p{\proj{2}}_*(\proj{1}^*\mathcal{E}\otimes\sheaf{P})
\otimes\kappa(y) \to
\Hypcoh{p}(X_y,(\proj{1}^*\sheaf{E}\otimes\sheaf{P})_{y})
$$
is trivially surjective, and by the base change theorem in fact
isomorphic, for all $p \neq n$. As the hyper direct images of a
complex of coherent sheaves are coherent for a proper map, we have
$$
\mathbf{R}^p{\proj{2}}_*(\proj{1}^*\mathcal{E}\otimes\sheaf{P})=0
$$
for $p\neq n$ by Nakayama's lemma. This proves the first part of
the proposition.

Now in particular
$\mathbf{R}^{n+1}{\proj{2}}_*(\proj{1}^*\mathcal{E}\otimes\sheaf{P})=0$.
Thus, by the second part of the base change theorem,
$\varphi^{n}(y)$ is an isomorphism. But as $\varphi^{n-1}(y)$ is
also surjective and thus isomorphic,
$\mathbf{R}^{n}{\proj{2}}_*(\proj{1}^*\mathcal{E}\otimes\sheaf{P})$
is free in a neighbourhood of $y$, again by the second part of
 \Ref{lem:basechange}.
\end{proof}

\begin{pro}\label{pro:basech&IT}
Let $X$, $Y$ and $S$ be as in \Ref{notation:WIT}, and let $u\maps
T \to S$ be a morphism of schemes. Suppose  that $\sheaf{E}$ is an
$IT(n)$-complex on $X$. Then, in the notation of
\Ref{pro:ITF-basechange}, $\Left i_X^*\sheaf{E}$ is a
$WIT(n)$-complex with respect to the pull-back $j^*\sheaf{P}$ of
$\sheaf{P}$ to $(X\times_S Y)_{(T)}$. Furthermore, if
$\widehat{\Left i_X^*\sheaf{E}}$ denotes the corresponding Fourier
transform, then
$$
i_Y^*\left(\widehat{\sheaf{E}}\right) = \widehat{\Left
i_X^*\sheaf{E}}.
$$
\end{pro}

\begin{proof}
By the assumptions and \Ref{pro:IT}, $\ITF{\sheaf{P}}{X\to
Y/S}(\sheaf{E})$ is a locally free sheaf shifted $n$ places to
the right. Hence \Ref{pro:ITF-basechange} gives
$$
i_Y^*\left(\ITF{\sheaf{P}}{X\to Y/S}(\sheaf{E})\right) =
\ITF{j^*\sheaf{P}}{X_{(T)}\to Y_{(T)}/T}(\Left i_X^*\sheaf{E}).
$$
But this shows that $\ITF{j^*\sheaf{P}}{X_{(T)}\to
Y_{(T)}/T}(\Left i_X^*\sheaf{E})$ is also a locally free sheaf
shifted $n$ places to the right. Both statements of the
proposition are now immediate.
\end{proof}

\subsection{Fourier transform for curves}\label{section:Fourier}
To fix terminology and notation, we first recall some basic facts
about Jacobians of curves; for details, see Milne
\cite{milne:abelian, milne:jacobian}.

\begin{notation}
Let $X$ be a smooth projective curve  of genus $g$. We denote by
$\Jac{X}$ a {\em Jacobian} of $X$, i.e., a scheme representing
the functor $T\mapsto\Pico{X/T}$. Let $\sheaf{M}$ be the
corresponding {\em universal sheaf} on $X\times\Jac{X}$. Recall
that $\Jac{X}$ is an Abelian variety of dimension $g$; let
$\widehat{\Jac{X}}$ denote its dual Abelian variety, and let
$\sheaf{P}$ be the {\em Poincar\'{e} sheaf} on
$\Jac{X}\times\widehat{\Jac{X}}$, normalised as in
\Ref{exa:relative_mukai}.
\end{notation}
\begin{para}\label{para:basepoint}
Choosing a base point $P\in X$ gives  the Abel-Jacobi  map
$i_{\!P}\maps X\to\Jac{X}$, taking the base point to $0$. Notice
that $i_{\!P}$ is a closed immersion. Furthermore, this choice
gives $\Jac{X}$ a principal polarisation and hence an isomorphism
$\varphi_P\maps\Jac{X}\xrightarrow{\sim}\widehat{\Jac{X}}$, which
we use henceforth to identify $\Jac{X}$ with its dual. Under this
identification, the pull-back
$(i_P\times\id{\Jac{X}})^*\sheaf{P}$ is just the universal sheaf
$\sheaf{M}$ on $X\times\Jac{X}$.
\end{para}
\begin{para}
Let $S$ be a separated $k$-scheme, $X_S=X\times S$, and let
$\Jac{X}_S=\Jac{X}\times S$ be the relative Jacobian of the
trivial family $X_S$. We have a Cartesian square
\begin{equation*}
\begin{CD}
X\times \Jac{X} \times S        @>{\proj{2}}>>        \Jac{X}_S \\
@V{\proj{1}}VV                                            @VV{}V    \\
X_S                     @>>{}>        S.    \\
\end{CD}
\end{equation*}
Let $\sheaf{M}_S$ be the pull-back of $\sheaf{M}$ to $X\times
\Jac{X} \times S$. The relative integral transform functor
$\ITF{\sheaf{M}_S}{X_S\to\Jac{X}_S/S}\maps \D^b_{coh}(X_S)\to
\D^b_{coh}(\Jac{X}_S)$ is given by
$$
    \ITF{\sheaf{M}_S}{X_S\to\Jac{X}_S/S}(\bullet)=
    \R {\proj{2}}_*(\proj{1}^*(\bullet)\otimes \sheaf{M}_S),
$$
where we can use the ordinary tensor product since $\sheaf{M}_S$
is locally free.
\end{para}
\begin{defn}\label{defn:Fourier}
The relative integral transform
$\ITF{\sheaf{M}_S}{X_S\to\Jac{X}_S/S}$ is called the {\em
relative Fourier functor} on $X\times S$ and is denoted by
$\Fourier{S}$. If $\sheaf{E}$ is $WIT$ with respect to
$\Fourier{S}$, the integral transform $\widehat{\sheaf{E}}$ is
called the {\em Fourier transform} of $\sheaf{E}$.
\end{defn}
\begin{pro}\label{pro:F->M}
Let $\Mukai{S}\maps\D^b_{coh}(\Jac{X}\times S)\to
\D^b_{coh}(\Jac{X}\times S)$ denote the relative Mukai transform.
Then
$$
\Fourier{S}=\Mukai{S} \circ (i_{P}\times\id{S})_*.
$$
\end{pro}
\begin{proof}
Consider the diagram
\begin{equation*}
\begin{CD}
X_S\times_S\Jac{X}_S        @>{j}>> \Jac{X}_S\times_S\Jac{X}_S
@>{p_2}>> \Jac{X}_S
\\
  @V{\proj{1}}VV @V{p_1}VV @VVV
\\
 X_S     @>>{i_P\times\id{S}}>   \Jac{X}_S   @>>> S,
\end{CD}
\end{equation*}
where the right-hand square is the fibre-product diagram and
$j=(i_P\times\id{S})\times_S\id{\Jac{X}_S}$. It is clear that the
left-hand square is also commutative, and that the composition of
the two top arrows is just the canonical projection $\proj{2}$.
But this means that the big rectangle is Cartesian, and hence so
is the left-hand square too.

By definition,
$$
\Mukai{S}(\bullet)= \R
{p_2}_*\left(p_1^*(\bullet)\otimes\sheaf{P}_S\right),
$$
where $\sheaf{P}_S$ is the pull-back of the Poincar\'{e} sheaf
onto $\Jac{X}_S\times_S\Jac{X}_S$. Clearly
$\sheaf{M}_S=j^*\sheaf{P}_S$.  Now by the projection formula
$$
\R j_* (\bullet\otimes\sheaf{M}_S) = \R j_*(\bullet) \otimes
\sheaf{P}_S.
$$

Because $p_1$ is flat as a base extension of a flat morphism, we
can do a base change \Ref{cor:basechange} around the left-hand
square to get
$$
p_1^*\circ\R (i_P\times\id{S})_* = \R j_*\circ \proj{1}^*.
$$
But $i_P\times\id{S}$ is a closed immersion and thus $\R
(i_P\times\id{S})_*=(i_P\times\id{S})_*$. Putting these
observations together, we get
\begin{align*}
\Mukai{S}((i_P\times\id{S})_*(\bullet))
    &= \R {p_2}_* \left(p_1^*
    \left((i_P\times\id{S})_*(\bullet)\right)
        \otimes\sheaf{P}_S\right) \\
    &= \R {p_2}_* \left( \R j_*\left(\proj{1}^*(\bullet)\right)
        \otimes\sheaf{P}_S\right) \\
    &= \R {p_2}_* \left(\R j_*
        \left(\proj{1}^*(\bullet)\otimes\sheaf{M}_S\right)\right) \\
    &= \R {\proj{2}}_* \left(\proj{1}^*(\bullet)\otimes\sheaf{M}_S\right)
    = \Fourier{S}(\bullet).
\end{align*}
\end{proof}

\begin{pro}\label{pro:hypcoh_of_fourier}
Let $X$ be a curve of genus $g$ and choose a base point $P\in X$
as in \Ref{para:basepoint}; we suppose made the identifications
given {\em loc. cit.} Let $S$ be a   $k$-scheme, and denote by
$j$ the embedding $S\cong(X\times S)_P\to X\times S$ of the fibre
over $P$. Let $\complex{\sheaf{E}}$ be a bounded complex of
locally free sheaves on $X\times S$. Then
$$
\Hypcoh{p}(\Jac{X}\times S, \Fourier{S}(\complex{\sheaf{E}})) =
\bigoplus_{i=1}^g
\Hypcoh{p-i}(S_P,j^*\complex{\sheaf{E}})^{\oplus\textstyle{\scriptstyle\binom{g-1}{i-1}}}.
$$
\end{pro}
\begin{proof}
By \Ref{pro:RGamma}   we have natural isomorphisms
$$
\Hypcoh{p}(\Jac{X}\times S, \Fourier{S}(\complex{\sheaf{E}})) =
\Hypcoh{p}(X\times S,
\complex{\sheaf{E}}\otimes\R\proj{1*}\sheaf{M}_{(S)})
$$
for all $p$.
\renewcommand{\theinlem}{\thethm.1}
\begin{inlem}\label{inlem:pr*M(s)}
With the notation of the proposition,
$\R\proj{1*}\sheaf{M}_{(S)}$ is the zero-differential complex
$\complex{\sheaf{C}}$ where $\sheaf{C}^i$ is the direct sum of
$\binom{g-1}{i-1}$ copies of $j_*\Ox{S}$ for $1\leq i \leq g$,
zero otherwise.
\end{inlem}

{\renewcommand{\theequation}{\theinlem.1}
Consider the Cartesian square
$$
\begin{CD}
X\times S \times \Jac{X}    @>{p'}>>    X\times\Jac{X}\\
@V{\proj{1}}VV                          @VVqV         \\
X\times S                   @>>{p}>     X.
\end{CD}
$$
By flat base change  around the square we get
\begin{align}\label{eq:lemma_on_Rpr*M}
\R\proj{1*}\sheaf{M}_{(S)}
    &=\R\proj{1*}{p'}^*\sheaf{M}
    =p^*\R q_*\sheaf{M}.
\end{align}
In order to compute $\R q_*\sheaf{M}$ on $X$, we consider the
Cartesian square
$$
\begin{CD}
X\times\Jac{X}  @>{i_P\times\id{}}>>      \Jac{X}\times\Jac{X}\\
@VqVV                       @VV{\pi_1}V                \\
X               @>>{i_P}>   \Jac{X}.
\end{CD}
$$
Now by the general base-change \Ref{cor:basechange} we have
\begin{align*}
\R q_*\sheaf{M}     &=\R q_*(i_P\times\id{})^*\sheaf{P}\\
                    &=\Left i_{\!P}^*\,\R\pi_{1*}\sheaf{P}.
\end{align*}
But $\R\pi_{1*}\sheaf{P}=k(0)[-g]$, the skyscraper sheaf at $0$
shifted $g$ places to the right (see the proof of the theorem of
\S13 in Mumford \cite{mumford:abelian}). Notice that $i_{\!P}$ is
a regular embedding; using Koszul resolutions it follows that
$\Left
i_{\!P}^*i_{P*}\Ox{X}=\bigwedge^{\bullet}\sheaf{N}_{X/\Jac{X}}$,
the zero-differential exterior-algebra complex of the conormal
sheaf of $X$ in $\Jac{X}$, concentrated in degrees $-g+1$ to $0$.
Similarly $\Left
i_{\!P}^*\,k(0)=\bigwedge^{\bullet}\sheaf{N}_{X/\Jac{X}}(P)$, the
exterior algebra of the fibre at $P$, whence the lemma follows
immediately taking into account the shift by $-g$. }

\vspace{5pt} Using the  projection formula  we have
$$
\Hypcoh{p}(X\times S,\complex{\sheaf{E}}\otimes j_*\Ox{S_P}) =
\Hypcoh{p}(S_P,j^*\complex{\sheaf{E}}).
$$
The proposition now follows from the lemma because
hypercohomology commutes with  direct sums.
\end{proof}

\section{Transforms of Higgs bundles}
We shall now apply the Fourier-transform machinery developed in
the previous section to stable Higgs bundles on curves.

\subsection{Definitions and basic properties}
\begin{defn}\label{defn:higgs}
A {\em Higgs bundle} on a smooth projective curve is a pair
$\Higgs{E}=(\sheaf{E},\theta)$, where $\sheaf{E}$ is a locally free sheaf on
$X$, and $\theta$ is a morphism
$\sheaf{E}\to\sheaf{E}\otimes\omega_X$. The morphism $\theta$ is
often called the {\em Higgs field}.
The Higgs bundle $\Ox{X}\xrightarrow{0}\omega_X$ is called
{\em trivial}.

The rank and degree (i.e., the first Chern class) of a Higgs
bundle $(\sheaf{E},\theta)$ mean the rank and degree of the
underlying sheaf $\sheaf{E}$. If $\Higgs{E}=\higgs{E}{\theta}$ and
$\Higgs{F}=\higgs{F}{\eta}$ are Higgs bundles, by a {\em
morphism} $\Higgs{E}\to\Higgs{F}$ we understand a morphism of
sheaves $\varphi\maps\sheaf{E}\to\sheaf{F}$ making the square
$$
\begin{CD}
\sheaf{E}       @>{\theta}>> \sheaf{E}\otimes\omega_X\\
@V{\varphi}VV                @VV{\varphi\otimes\id{}}V\\
\sheaf{F}       @>>{\eta}>   \sheaf{F}\otimes\omega_X
\end{CD}
$$
commutative.
\end{defn}

\begin{para}
Let $\Higgs{E}=\higgs{E}{\theta}$ be a Higgs bundle on $X$. Then
we can consider it as a complex of sheaves concentrated in
degrees $0$ and $1$, and hence as an object in $\D_{coh}^b(X)$.
When we write $\Higgs{E}\otimes \sheaf{F}$ or
$\Hypcoh{\bullet}(X,\Higgs{E})$ etc., we consider the Higgs
bundle as a sheaf complex this way. Notice that the image of
$\Higgs{E}$ in $\D_{coh}^b(X)$ does not uniquely determine the
isomorphism class of the {\em Higgs bundle} $\higgs{E}{\theta}$.
In fact, multiplying $\theta$ by a non-zero constant gives a
quasi-isomorphic complex; however, the resulting Higgs bundle is
not in general isomorphic.
\end{para}

\begin{defn}\label{defn:stability}
A Higgs bundle $\higgs{E}{\theta}$ is called {\em stable} if for
any locally free subsheaf $\sheaf{F}$ of $\sheaf{E}$ satisfying
$\theta(\sheaf{F})\subset\sheaf{F}\otimes\omega_X$,
we have
$$
\frac{\deg{\sheaf{F}}}{\rank{\sheaf{F}}} <
\frac{\deg{\sheaf{E}}}{\rank{\sheaf{E}}}.
$$
\end{defn}

\begin{thm}\label{thm:hausel-vanishing}
Let $\Higgs{E}=\higgs{E}{\theta}$ be a non-trivial stable Higgs bundle on $X$
with $\deg(\Higgs{E})=0$. Then
$$
    \Hypcoh{p}(X,\Higgs{E})=0
$$
for $p\neq 1$.
\end{thm}
\begin{proof}
Hausel \cite{hausel:vanishing} Corollary (5.1.4.). Notice that
$\Hypcoh{p}(X,\Higgs{E})=0$ automatically for  $p > 2$ because
$\dim(X)=1$ and the length of the complex $\Higgs{E}$ is $2$.
\end{proof}

\begin{pro}\label{pro:stability&Pico}
If a Higgs bundle $\Higgs{E}$ is stable, then so is
$\Higgs{E}\otimes \sheaf{L}$, where $\sheaf{L}$ is an element of
$\Pico{X}$.
\end{pro}
\begin{proof}
Let $\sheaf{F}\subset \sheaf{E}\otimes\sheaf{L}$ be a subbundle
stable under $\theta\otimes{\id{}}_{\sheaf{L}}$. Then
$\sheaf{F}\otimes\sheaf{L}^{-1}$ is a subbundle of $\sheaf{E}$
stable under $\theta$. But tensoring with $\sheaf{L}$ affects
neither the ranks nor the degrees of $\sheaf{E}$ and $\sheaf{F}$,
and hence the lemma follows from the stability of the Higgs
bundle $\Higgs{E}$.
\end{proof}

\begin{para}\label{para:E(a)}
Let $\Higgs{E}=\higgs{E}{\theta}$ be a Higgs bundle and
$\alpha\in H^0(X,\omega_X)$ a global $1$-form. Then
$\id{\sheaf{E}}\otimes\alpha$ is canonically identified with a
morphism $\sheaf{E}\to\sheaf{E}\otimes\omega_X$. We denote the
Higgs bundle $\higgs{E}{\theta+\id{\sheaf{E}}\otimes\alpha}$ by
$\Higgs{E}(\alpha)$.
\end{para}

\begin{lem}\label{lem:stability&forms}
Let $\Higgs{E}$ be a stable Higgs bundle. Then $\Higgs{E}(\alpha)$
is also stable for any $\alpha\in H^0(X,\omega_X)$.
\end{lem}
\begin{proof}
Let $\sheaf{F}\subset \sheaf{E}$ be a subbundle stable under
$\theta_{\alpha}=\theta+\id{}\otimes\alpha$. Let $t\in\Gamma(U,\sheaf{F})$.
Then
$\theta_{\alpha}(t)=\theta(t)+t\otimes\alpha\in\Gamma(U,\sheaf{F}\otimes\omega_X)$.
But $t\otimes\alpha\in\Gamma(U,\sheaf{F}\otimes\omega_X)$ too, and hence
$\theta(t)\in \Gamma(U,\sheaf{F}\otimes\omega_X)$. Thus $\sheaf{F}$ is stable
under $\theta$, and the lemma follows from the stability of
$\Higgs{E}$.
\end{proof}
We shall now introduce an important construction of algebraic
families of Higgs bundles. For details about projective bundles
see for example EGA II \cite{EGAII} \S4.
\begin{para}\label{para:THC}
Let $\Higgs{E}=\higgs{E}{\theta}$ be a Higgs bundle on a curve
$X$ of genus $g$, and let $\pi\maps X\to\Spec{k}$ be the
structural morphism. Then the $k$-rational points of the vector
bundle (or affine space) $\mathbf{V}((\pi_*\omega_X)^{\vee})$ are
canonically identified with the elements of $H^0(X,\omega_X)$; we
use the notation $H^0(X,\omega_X)$ also for this scheme if no
confusion seems likely. Let
$\sheaf{D}=\pi^*((\pi_*\omega_X)^{\vee})=(\pi^*\pi_*\omega_X)^{\vee}$;
we have the canonical adjunction morphism
$$
\varphi\maps\sheaf{D}^{\vee}=\pi^*\pi_*\omega_X\to\omega_X.
$$
Let
$\tilde{\varphi}\maps\sheaf{D}^{\vee}\to\omega_X\otimes\sheaf{E}nd(\sheaf{E})$
be the morphism
$$
t\mapsto\varphi(t)\otimes\id{\sheaf{E}}.
$$
On the other hand, let
$\psi\maps\Ox{X}\to\omega_X\otimes\sheaf{E}nd(\sheaf{E})$
be the map that takes $1$ to $\theta$. Putting these together we
get a morphism
$$
\gamma=\tilde{\varphi}+\psi\maps\sheaf{D}^{\vee}\oplus\Ox{X}\to\omega_X\otimes\sheaf{E}nd(\sheaf{E}).
$$
Because $\sheaf{D}\oplus\Ox{X}=\pi^*((\pi_*\omega_X)^{\vee}\oplus
k)$, we have a canonical isomorphism
$$
\Proj[X]{}(\sheaf{D}\oplus\Ox{X})=X\times
\Proj[k]{}((\pi_*\omega_X)^{\vee}\oplus k)
=X\times\Proj{}(H^0(X,\omega_X)\oplus k) \cong X\times\Proj[k]{g}.
$$
Let $p\maps P=\Proj[X]{}(\sheaf{D}\oplus\Ox{X})\to X$ be the
projection. There is  the canonical surjection
$p^*(\sheaf{D}\oplus\Ox{X})\to\Ox{P}(1)$, and so  by
dualising a canonical
$\Ox{P}(-1)\to p^*(\sheaf{D}^{\vee}\oplus\Ox{X})$. Composing this
morphism with $p^*\gamma$ we get a morphism
$$
\Ox{P}(-1)\to p^*(\omega_X\otimes\sheaf{E}nd(\sheaf{E})),
$$
or in other words a global section of
$p^*(\omega_X\otimes\sheaf{E}nd(\sheaf{E}))\otimes\Ox{P}(1)$. We
interpret this section as a morphism
$$
\Theta\maps p^*\sheaf{E}\to p^*\sheaf{E}\otimes
p^*\omega_X\otimes\Ox{P}(1),
$$
and denote this complex of sheaves (in degrees $0$ and $1$) on $P$
by $\THC{E}$.

In more pedestrian terms, let $(\alpha_i)_i$ be a basis of
$H^0(X,\omega_X)$, and let $(\alpha_i^*)_i$ be the dual basis of
$H^0(X,\omega_X)^{\vee}$. Let $t\maps k\to k$ be the canonical
coordinate on $k$; then $(t,\alpha_1^*,\dotsc,\alpha_g^*)$ forms a
basis of the global sections of $\Ox{\Proj{g}}(1)$, and
$H^0(X,\omega_X)$ corresponds to the open affine subscheme of
$\Proj{g}$ with $t\neq 0$. Now
$$
\Theta= \theta\otimes t+\sum_{i=1}^g
\alpha_i\otimes\id{}\otimes\alpha_i^*.
$$
\end{para}

\begin{rem}\label{rem:E(a)}
Notice that for $\alpha\in H^0(X,\omega_X)$ the restriction of
$\THC{E}$ to $X\times\set{\alpha}$ is just $\Higgs{E}(\alpha)$ of
\Ref{para:E(a)}.
\end{rem}

\begin{pro}\label{pro:TFT}
Let $\Higgs{E}$ be a  stable Higgs bundle of degree $0$ and rank
$\geq 2$ on a curve $X$ of genus $g\geq 2$. Then the complex
$\THC{E}$ on $X\times\Proj{g}$ is $WIT(1)$ with respect to the
relative Fourier functor
$\Fourier{\Proj{g}}\maps\D^b_{coh}(X\times\Proj{g})\to\D^b_{coh}(\Jac{X}\times\Proj{g})$.
Moreover, the Fourier transform $\left(\THC{E}\right)
\widehat{\big{.}}$ is a locally free sheaf on
$\Jac{X}\times\Proj{g}$.
\end{pro}
\begin{proof}
By \Ref{pro:IT} we are reduced to showing that $\THC{E}$ is
$IT(1)$ with respect to $\sheaf{M}_{(\Proj{g})}$. We  consider
two cases. Let $U$ denote the open subset $H^0(X,\omega_X)$ in
$\Proj{g}$.

\vspace{5pt} A) Let $(\xi,\alpha)\in\Jac{X}\times U$. Then (using
the notation of \Ref{para:Fs(k(x))})
$$
\left(\THC{E}\right)_{(\xi,\alpha)} \cong \Higgs{E}(\alpha),
$$
and we need to show that
$$
\Hypcoh{p}(X,\Higgs{E}(\alpha)\otimes\sheaf{M}_{\xi})=0
$$
for $p\neq 1$. But this follows from \Ref{pro:stability&Pico},
\Ref{lem:stability&forms} and \Ref{thm:hausel-vanishing}. Notice
that for a rank-$1$ Higgs bundle $\Higgs{E}$ one of the bundles
$\Higgs{E}(\alpha)$ would be trivial, and the vanishing theorem
\Ref{thm:hausel-vanishing} would fail.

\vspace{5pt} B) Let $(\xi,z)\in\Jac{X}\times(\Proj{g}-U)$. We
consider the second hypercohomology spectral sequence:
$$
    \SSeqII_2^{pq} = H^p(X,H^q(\left(\THC{E}\right)_{(\xi,z)}\otimes\sheaf{M}_{\xi}))
    \Rightarrow\Hypcoh{p+q}(X,\left(\THC{E}\right)_{(\xi,z)}\otimes\sheaf{M}_{\xi}).
$$
But
$$
\left(\THC{E}\right)_{(\xi,z)} \cong \higgs{E}{\id{}\otimes\alpha}
$$
for a $1$-form $\alpha\neq 0$, determined up to multiplication by
a non-zero scalar . Now $\id{}\otimes\alpha$ is clearly an
injective map of sheaves; let $\sheaf{S}$ be its cokernel. Thus
the $E_2$-terms of the spectral sequence are
$$
\specseq{\SSeqII_2^{pq}}{H^0(X,\sheaf{S}\otimes\sheaf{M}_{\xi})}
{H^1(X,\sheaf{S}\otimes\sheaf{M}_{\xi})}{0}{0}{}{}{}{}
$$
But $\sheaf{S}$ is a direct sum of skyscraper sheaves supported on
the divisor of zeroes of the one-form, and since skyscraper
sheaves are flasque, we have
$H^1(X,\sheaf{S}\otimes\sheaf{M}_{\xi})=0$. Hence
$\Hypcoh{0}(X,\left(\THC{E}\right)_{(\xi,z)}\otimes\sheaf{M}_{\xi})
=\Hypcoh{2}(X,\left(\THC{E}\right)_{(\xi,z)}\otimes\sheaf{M}_{\xi})
= 0$.
\end{proof}

\begin{defn}
Let $\Higgs{E}$ be a  stable Higgs bundle of degree $0$ and rank
$r\geq 2$ on a curve $X$ of genus $g\geq 2$. Then the locally free
sheaf $\left(\THC{E}\right) \widehat{\big{.}}$ on
$\Jac{X}\times\Proj[k]{g}$ is called (by abuse of language) the
{\em total Fourier transform} of $\Higgs{E}$ and is denoted by
$\TFT{E}$.
\end{defn}

\begin{pro}\label{pro:FT_of_Higgs}
Let $\Higgs{E}$ and $X$ be as in \Ref{pro:TFT}, and let
$\alpha\in H^0(X,\omega_X)$. Then
$$
\TFT{\Higgs{E}}_{\alpha}\cong\widehat{\Higgs{E}(\alpha)},
$$
where the left-hand side denotes the absolute Fourier transform.
\end{pro}
\begin{proof}
By the proof of \Ref{pro:TFT} $\THC{E}$ is $IT(1)$. Now the
proposition follows from Remark \Ref{rem:E(a)} and Proposition
\Ref{pro:basech&IT} applied to the immersion $\set{\alpha}\to
H^0(X,\omega_X)\to\Proj{g}$.
\end{proof}

\begin{pro}\label{pro:rank}
Let $\Higgs{E}=\higgs{E}{\theta}$ be a non-trivial stable Higgs
bundle of degree $0$ on a curve $X$ of genus $g\geq 2$. Then the
rank of the total Fourier transform $\TFT{\Higgs{E}}$ is $(2g-2)
\rank(\sheaf{E})$.
\end{pro}
\begin{proof}
It follows from \Ref{pro:TFT} and \Ref{pro:FT_of_Higgs} that
$\rank(\TFT{\Higgs{E}}) = \dim \Hypcoh{1}(X,\Higgs{E})$. Consider
the first hypercohomology spectral sequence
$$
    \SSeqI_2^{pq}= H^p(H^q(X,\Higgs{E})) \Rightarrow
    \Hypcoh{p+q}(X,\Higgs{E}).
$$
The $E_1$-terms of the sequence are:
$$
\specseq{\SSeqI_1^{pq}}{H^1(X,\sheaf{E})}{H^1(X,\sheaf{E}\otimes\omega_X)}
{H^0(X,\sheaf{E})}{H^0(X,\sheaf{E}\otimes\omega_X)}
{->}{H^1(\theta)}{->}{H^0(\theta)}
$$
The sequence clearly degenerates at $E_2$, i.e.,
$\SSeqI_{\infty}^{pq}=\SSeqI_2^{pq}$, and hence
\begin{align*}
    \SSeqI_2^{0,0}&\cong\Hypcoh{0}(X,\Higgs{E})\quad\text{and} \\
    \SSeqI_2^{1,1}&\cong\Hypcoh{2}(X,\Higgs{E}).
\end{align*}
But these hypercohomologies  vanish by
\Ref{thm:hausel-vanishing}, and thus $H^0(X,\theta)$ is injective
and $H^1(X,\theta)$ is surjective. On the other hand,
\begin{equation*}
\begin{split}
    \Hypcoh{1}(X,\Higgs{E})
       \cong\SSeqI_{\infty}^{0,1}\oplus\SSeqI_{\infty}^{1,0}
       = \ker{H^1(X,\theta)}\oplus \coker{H^0(X,\theta)},
\end{split}
\end{equation*}
and hence
\begin{equation*}
\begin{split}
\dim \Hypcoh{1}(X,\Higgs{E})
    &= \dim H^1(X,\sheaf{E})-\dim H^1(X,\sheaf{E}\otimes\omega_X)\\
    & \quad +\dim
    H^0(X,\sheaf{E}\otimes\omega_X)-\dim H^0(X,\sheaf{E}) \\
    &= \chi(\sheaf{E}\otimes\omega_X)-\chi(\sheaf{E}).
\end{split}
\end{equation*}
But as $\deg(\sheaf{E})=0$, the Riemann-Roch theorem gives
\begin{align*}
\chi(\sheaf{E})                & = (1-g)\rank(\sheaf{E})\quad\text{and}\\
\chi(\sheaf{E}\otimes\omega_X) & = (g-1)\rank(\sheaf{E}),
\end{align*}
whence the result follows immediately.
\end{proof}

\begin{pro}\label{pro:cohomology_of_Higgs}
Let $\Higgs{E}$ be a stable Higgs bundle of rank $r\geq 2$ and
degree $0$ on a curve $X$ of genus $g\geq 2$. Then
$$
\dim_k H^p(\Jac{X}\times\Proj{g},\TFT{\Higgs{E}})
=rg{\binom{g-1}{p-1}},
$$
when $1\leq p\leq g$, and zero otherwise.

\end{pro}
\begin{proof}
Let $P\in X$ be a base point giving an embedding $i_P\maps
X\to\Jac{X}$, and denote by $j$ the embedding $\Proj{g}\to
X\times\Proj{g}$ of the fibre $\proj{X}^{-1}(P)$. Then by
\Ref{pro:hypcoh_of_fourier}
\begin{equation}\label{eq:hypcoheq}
\begin{split}
H^p(\Jac{X}\times\Proj{g},\TFT{\Higgs{E}})
    &=\Hypcoh{p+1}(\Jac{X}\times\Proj{g},\Fourier{\Proj{}}(\THC{E}))\\
    &=\bigoplus_{i=1}^g\Hypcoh{p+1-i}(\Proj{g},j^*\THC{E})^{\oplus\textstyle{\scriptstyle\binom{g-1}{i-1}}}.
\end{split}
\end{equation}
We apply the first hypercohomology spectral sequence
$$
\SSeqI_2^{pq}=H^p(H^q(\Proj{g},j^*\THC{E}))\Rightarrow
\Hypcoh{p+q}(\Proj{g},j^*\THC{E}).
$$
The $E_1$-terms are given by
$$
\specseq{\SSeqI_1^{pq}}{H^1(\Proj{g},\Ox{\Proj{g}}^r)}{H^1(\Proj{g},{\Ox{\Proj{g}}(1)}^r)}
{H^0(\Proj{g},\Ox{\Proj{g}}^r)}{H^0(\Proj{g},{\Ox{\Proj{g}}(1)}^r).}{->}{}{->}{d}
$$
The standard results on the cohomology of a projective space
(Hartshorne \cite{AG} III.5.1) show that the
$E_1^{0,1}=E_1^{1,1}=0$. Furthermore, it is clear from the
definition \Ref{para:THC} of $\THC{E}$ that
$d=H^0(\Proj{g},j^*\Theta)$ is an injection. Thus we see that
$$
\dim\Hypcoh{p}(\Proj{g},j^*\THC{E})
    =\begin{cases}
        rg & \text{if }$p=1$,\\
        0   & \text{otherwise}.
     \end{cases}
$$
Thus in the direct sum of \Ref{eq:hypcoheq} we have non-zero
cohomology only when $i=p$, and the result follows immediately.
\end{proof}

\begin{pro}\label{prop:ch}
Let $\Higgs{E}=\higgs{E}{\theta}$ be a stable non-trivial Higgs
bundle on a smooth projective curve $X$ of genus $g\geq 2$, with
$r=\rank(\sheaf{E})\geq 2$ and $\deg(\sheaf{E})=0$. Then
$$
\ch(\TFT{\Higgs{E}})=\rank(\Higgs{E})
\left(g-1+(g-1)\proj{\Proj{}}^*\ch(\Ox{\Proj{g}}(1))+
t.(1-\proj{\Proj{}}^*\ch(\Ox{\Proj{g}}(1)))\right),
$$
where $t$ is the class of the $\Theta$-divisor on $\Jac{X}$.
\end{pro}

\begin{proof}
This is an easy application of the Grothendieck-Riemann-Roch
formula.
\end{proof}

\subsection{Invertibility}
\begin{thm}\label{thm:invertibility}
Let $\Higgs{E}$ and $\Higgs{F}$ be two Higgs bundles on a curve $X$ of genus $g\geq 2$.
If $\TFT{\Higgs{E}}\cong\TFT{\Higgs{F}}$, then
$\Higgs{E}\cong\Higgs{F}$ {\em as Higgs bundles}.
\end{thm}
\begin{proof}
We show this by actually exhibiting a process of recovering a
Higgs bundle $\Higgs{E}$ from its total Fourier transform
$\TFT{\Higgs{E}}$.

{\em Step $1$.}
Choose a base point $P\in X$ as in \Ref{para:basepoint}, and let
$i_P\maps X\to\Jac{X}$ be the corresponding embedding. Denote by
$j$ the immersion $i_P\times\id{\Jac{X}}$.
Then by \Ref{pro:F->M} $\Fourier{\Proj{g}}=\Mukai{\Proj{g}}\circ j_*$.
By \Ref{thm:relative_mukai_equivalence} $\Mukai{\Proj{g}}$ is a
category equivalence; let
$\mathbf{G}$
be its inverse.
Now by definition $\TFT{\Higgs{E}}=\Fourier{\Proj{g}}(\THC{E})[1]$,
and hence
$$
\mathbf{G}(\TFT{\Higgs{E}})[-1]= j_*(\THC{E}).
$$
\renewcommand{\theinlem}{\thethm.1}
\addtocounter{equation}{1}
\begin{inlem}\label{inlem:coker}
The differential $\Theta$ of the complex $\THC{E}$ is injective.
\end{inlem}
Let $U\subset X\times\Proj{g}$ be an open subset and
$s\in\Gamma(U,\proj{}^*\sheaf{E})$ a non-zero section. There is a
point $z=(x,p)\in U$ for which $s(z)\neq 0$. Because $\sheaf{E}$
is locally free, it follows (using Nakayama's lemma) that there is
an open neighbourhood $V\subset U$ of $z$ such that $s(z')\neq 0$
for $z'\in V$.  If $\Theta(z)(s(z))=0$, it follows from the
definition of $\Theta$ that there is a point $y\in V$ with
$\Theta(y)(s(y))\neq 0$, and  in particular $\Theta_U(s)\neq 0$.
But this shows that $\Theta$ is injective as a morphism of
presheaves and hence as a sheaf morphism too. Thus the lemma is
proved.

\vspace{5pt} By the lemma there is an exact sequence
\begin{equation}\label{eq:exseq}
0\to \proj{}^*\sheaf{E}
    \xrightarrow{\Theta}
    \proj{}^*(\sheaf{E}\otimes  \omega_X)
    \otimes
    \proj{}^*\sheaf{O}_{\Proj{g}}(1)
    \to \sheaf{R}\to 0,
\end{equation}
and consequently $\THC{\sheaf{E}}$ is quasi-isomorphic to
$\sheaf{R}[-1]$. It  follows from this that
$\mathbf{G}(\TFT{\Higgs{E}})= j_*\sheaf{R}$ in
$\D^b_{coh}(X\times\Proj{g})$. Since $j_*\sheaf{R}$ is an honest
sheaf, $\mathbf{G}(\TFT{\Higgs{E}})= j_*\sheaf{R}$ also in
$\Cat{Mod}(X\times\Proj{g})$. This means that we can recover the
cokernel $\sheaf{R}$ of $\THC{E}$ on $X\times\Proj{g}$ as
$j^*(\mathbf{G}(\TFT{\Higgs{E}}))$.

\vspace{5pt}
{\em Step $2$.}
Tensor \Ref{eq:exseq} with
$\proj{}^*\Ox{\Proj{g}}(-1)$ and obtain the exact sequence
\begin{equation}\label{eq:exseq2}
0\to \proj{}^*\sheaf{E}\otimes\proj{}^*\Ox{\Proj{g}}(-1)
    \xrightarrow{\Theta\otimes\id{}}
    \proj{}^*(\sheaf{E}\otimes \omega_X)
    \to \sheaf{R}\otimes\proj{}^*\Ox{\Proj{g}}(-1)
    \to 0.
\end{equation}
We shall use the long exact $\R\proj{X*}$-sequence associated to
\Ref{eq:exseq2}. By the projection formula
\begin{align*}
\R \proj{X*}(\proj{X}^*\sheaf{E}\otimes\proj{\Proj{}}^*\Ox{\Proj{g}}(-1))&=
\sheaf{E}\otimes\R\proj{X*}\proj{\Proj{}}^*\Ox{\Proj{}}(-1),
\quad\text{and}\\
\R \proj{X*}(\proj{X}^*(\sheaf{E}\otimes\omega_X))&=
\sheaf{E}\otimes\omega_X\otimes\R\proj{X*}\Ox{X\times\Proj{}}.
\end{align*}
Now it follows from base change and the standard formulas for the
cohomology of projective spaces  that
\begin{gather*}
\proj{X*} \proj{\Proj{}}^* \Ox{\Proj{}}(-1)=
R^1 \proj{X*} \proj{\Proj{}}^* \Ox{\Proj{}}(-1)=0,\quad\text{and}\\
\proj{X*}  \Ox{X\times\Proj{}}=\Ox{X}.
\end{gather*}
It follows then from the long exact sequence that
$\proj{X*}(\sheaf{R}\otimes\proj{\Proj{}}^*\Ox{\Proj{}}(-1))\cong\sheaf{E}\otimes\omega_X$,
and that  we may consequently recover the underlying sheaf
$\sheaf{E}$ of $\Higgs{E}$ from $\sheaf{R}$ by twisting by
$\Ox{\Proj{}}(-1)$, projecting down to $X$, and twisting by
$\omega_X^{\vee}$.

\vspace{5pt} {\em Step $3$.} It remains to recover the Higgs
field $\theta$. This will be done after discarding much of the
information contained in $\sheaf{R}$. We choose a non-zero
$\alpha\in H^0(X,\omega_X)$, and we let $U=\Spec{A}$ be an open
affine subscheme of $X$ over which $\alpha$ does not vanish; then
$\alpha$ gives a trivialisation of $\omega_X$ on $U$. Clearly it
is enough to recover $\theta$ over $U$.

Let $V$ be the subvectorspace of $H^0(X,\omega_X)$ generated by
$\alpha$. We can consider $V$ as a closed subscheme of the open
subscheme $H^0(X,\omega_X)$ of $\mathbf{P}(H^0(X,\omega_X)\oplus
k)$. Furthermore, we consider $U\times V$ as a subscheme of
$U\times \mathbf{P}(H^0(X,\omega_X)\oplus k)$, and let
$\sheaf{S}$ be the restriction of $\sheaf{R}$ to $U\times V$; it
is just the cokernel of $\Theta$ restricted to $U\times V$.
Notice that $U\times V\cong \Spec{A[T]}$.

On $U$ the underlying sheaf $\sheaf{E}$ of $\Higgs{E}$ corresponds
to an $A$-module $M$ and $\theta$ corresponds to an endomorphism
$u$ of $M$. Furthermore, the  pull-back of $\sheaf{E}$ to $U\times
V$ corresponds to $M[T]=M\otimes_A A[T]$. By the definition of
$\Theta$ \Ref{para:THC}, $\Theta|_{U\times V}$ corresponds to the
$A[T]$-linear map
$$
\psi=\id{M}\otimes T+u\otimes \id{A[T]}.
$$
But $\psi$ fits into the exact sequence
$$
M[T] \xrightarrow{\psi} M[T] \rightarrow M_u \rightarrow 0,
$$
where $M_u$ is the $A[T]$-module with $T$ acting on $M$ as $u$
(cf. Bourbaki \cite{ALG1-3}, Ch. III \S 8 no. 10). Hence
$\sheaf{S}=(M_u)^{\sim}$. But the structure of $A[T]$-module of
$M_u$ determines $u$ and hence $\theta|_U$.
\end{proof}

\begin{rem}
Lemma 6.8 in Simpson \cite{simpson:moduli2} gives a description of
Higgs bundles on $X$ as coherent sheaves on the total space of the
cotangent bundle of $X$. The scheme $U\times V$ in Step 3 of the
proof is the total space of the cotangent bundle of $U$, and the
coherent sheaf $\sheaf{S}$ on $U\times V$ is  the one that
corresponds to $\Higgs{E}|_{U}$ under Simpson's correspondence.
\end{rem}

\begin{cor} The functor $\mathbf{TF}$ from the category of stable non-trivial Higgs
bundles on $X$ with vanishing Chern classes to
$\Cat{Mod}(\Jac{X}\times\Proj{g})$ is fully faithful.
\end{cor}

\begin{proof}
Let $\Higgs{E}$ and $\Higgs{E}'$ be Higgs bundles on $X$ and let
$\sheaf{R}$ and $\sheaf{R}\,'$ be the cokernels of
$\THC{\Higgs{E}}$ and $\THC{\Higgs{E}'}$ respectively. Because the
relative Mukai transform is an equivalence of categories, we have
$$
\Hom(\TFT{\Higgs{E}},\TFT{\Higgs{E}'})=
\Hom(\sheaf{R},\sheaf{R}\,').
$$
Thus faithfulness is clear. On the other hand, let
$\varphi\maps\sheaf{R}\to\sheaf{R}'$; using the notation of the
proof of the theorem, the previous remark shows that
$\varphi|_{U\times V}$ gives a morphism of Higgs bundles
$\Higgs{E}|_U\to{\Higgs{E}'}|_U$. But as the genus of $X$ is at
least $2$, the canonical linear system $\abs{\omega_X}$ has no
base points. Hence we can cover $X$ by open sets like $U$; it is
clear that the morphisms thus obtained glue to give a morphism
$\Higgs{E}\to\Higgs{E}'$.
\end{proof}

\bibliographystyle{std}
\bibliography{maths}

\newcommand{\noopsort}[1]{}
\begin{thebibliography}{10}

\bibitem{bondal-orlov:semiorthogonal}
A.~Bondal and D.~Orlov: Semiorthogonal decomposition for algebraic varieties
  (1995), preprint alg-geom/9506012.

\bibitem{bondal-orlov:reconstructing}
--- Reconstruction of a variety from the derived category and groups of
  autoequivalences (1997), preprint alg-geom/9712029.

\bibitem{ALG1-3}
Nicolas Bourbaki: \emph{Alg\`ebre {Ch.} 1--3}, Hermann, Paris, 1970.

\bibitem{bridgeland:thesis}
Tom Bridgeland: \emph{Fourier-Mukai transforms for surfaces and moduli spaces
  of stable sheaves}, Ph.D. thesis, University of Edinburgh (1998).

\bibitem{bridgeland:equivalences}
--- Equivalences of triangulated categories and {F}ourier-{M}ukai transforms,
  \emph{Bull. London Math. Soc.} \textbf{31}, no.~1 (1999) 25--34.

\bibitem{Donaldson-Kronheimer}
Simon~K. Donaldson and Peter~B. Kronheimer: \emph{The Geometry of
  Four-Manifolds}, Oxford Mathematical Monographs, Clarendon Press, Oxford,
  1990.

\bibitem{gelfand-manin:methods}
Sergei Gelfand and Yuri Manin: \emph{Methods of Homological Algebra},
  Springer-Verlag, Berlin, 1996.

\bibitem{EGAII}
Alexander Grothendieck and Jean Dieudonn\'{e}: Él\'{e}m\'{e}nts de
  g\'{e}om\'{e}trie alg\'{e}brique {\noopsort{c}} {II}, \emph{Inst. Hautes
  \'Etudes Sci. Publ. Math.} \textbf{8} (\noopsort{ac}1961) $\negmedspace$.

\bibitem{EGAIII}
--- Él\'{e}m\'{e}nts de g\'{e}om\'{e}trie alg\'{e}brique {\noopsort{d}} {III},
  \emph{Inst. Hautes \'Etudes Sci. Publ. Math.} \textbf{{11,17}}
  (\noopsort{b}{1961,1963}) $\negmedspace$.

\bibitem{hartshorne:residues}
Robin Hartshorne: \emph{Residues and Duality}, Lecture Notes in Mathematics 20,
  Springer-Verlag, Berlin, 1966.

\bibitem{AG}
--- \emph{Algebraic Geometry}, Graduate Texts in Mathematics 52,
  Springer-Verlag, New York, 1977.

\bibitem{hausel:vanishing}
Tam\'{a}s Hausel: Vanishing of intersection numbers on the moduli space of
  {Higgs} bundles, \emph{Adv. Theor. Math. Phys.} \textbf{2}, no.~5 (1998)
  1011--1040.

\bibitem{hitchin:self-duality}
Nigel~J. Hitchin: Self-duality equations on a {Riemann} surface, \emph{Proc.
  London Math. Soc. (3)} \textbf{55} (1987) 59--126.

\bibitem{hitchin:dirac}
--- The {D}irac operator, in Martin Bridson and Simon Salamon (Editors),
  \emph{Topics in Geometry and Topology}, Oxford University Press, Oxford, to
  appear.

\bibitem{illusie:categories}
Luc Illusie: Cat\'egories d\'eriv\'ees et dualit\'e: travaux de {J}.-{L}.
  {V}erdier, \emph{Enseign. Math. (2)} \textbf{36}, no. 3-4 (1990) 369--391.

\bibitem{jardim:thesis}
Marcos Jardim: \emph{Nahm Transform for Doubly Periodic Instantons}, Ph.D.
  thesis, University of Oxford (1999), available as math.DG/9912028.

\bibitem{jardim:constructing}
--- Construction of doubly-periodic instantons, \emph{Comm. Math. Phys.}
  \textbf{216}, no.~1 (2001) 1--15.

\bibitem{kashiwara-shapira:sheaves}
Masaki Kashiwara and Pierre Schapira: \emph{Sheaves on Manifolds}, Grundlehren
  292, Springer-Verlag, Berlin, 1990.

\bibitem{maciocia:generalized}
Antony Maciocia: Generalized {F}ourier-{M}ukai transforms, \emph{J. Reine
  Angew. Math.} \textbf{480} (1996) 197--211.

\bibitem{milne:abelian}
John~S. Milne: Abelian varieties, in Gary Cornell and Joseph~H. Silverman
  (Editors), \emph{Arithmetic Geometry (Storrs, 1984)}, Springer-Verlag, New
  York, 1986, (pp. 103--150).

\bibitem{milne:jacobian}
--- Jacobian varieties, in Gary Cornell and Joseph~H. Silverman (Editors),
  \emph{Arithmetic Geometry (Storrs, 1984)}, Springer-Verlag, New York, 1986,
  (pp. 167--212).

\bibitem{mukai:duality}
Shigeru Mukai: Duality between {$\mathbf{D}(X)$} and $\mathbf{D}(\hat{X})$ with
  its application to {Picard} sheaves, \emph{Nagoya Math. J.} \textbf{84}
  (1984) 153--175.

\bibitem{mukai:fourier}
--- Fourier functor and its application to the moduli of bundles on an abelian
  variety, in \emph{Algebraic geometry, Sendai, 1985}, North-Holland,
  Amsterdam, 1987, (pp. 515--550).

\bibitem{mumford:abelian}
David Mumford: \emph{Abelian Varieties}, Oxford University Press, Oxford, 1970.

\bibitem{GAGA}
Jean-Pierre Serre: G\'{e}om\'{e}trie alg\'{e}brique et g\'{e}om\'{e}trie
  analytique, \emph{Ann. Inst. {Fourier}} \textbf{6} (1956) 1--42.

\bibitem{simpson:constructing}
Carlos~T. Simpson: Constructing variations of {H}odge structure using
  {Y}ang-{M}ills theory and applications to uniformization, \emph{J. Amer.
  Math. Soc.} \textbf{1}, no.~4 (1988) 867--918.

\bibitem{simpson:moduli2}
--- Moduli of representations of the fundamental group of a smooth projective
  variety. {I}{I}, \emph{Inst. Hautes \'Etudes Sci. Publ. Math.} \textbf{80}
  (1994) 5--79 (1995).

\bibitem{simpson:hodge}
--- The {H}odge filtration on nonabelian cohomology, in \emph{Algebraic
  geometry---Santa Cruz 1995}, Amer. Math. Soc., Providence, RI, 1997, (pp.
  217--281).

\bibitem{verdier:thesis}
Jean-Louis Verdier: Des cat\'{e}gories d\'{e}riv\'{e}es des cat\'{e}gories
  ab\'{e}liennes, \emph{Ast\`{e}risque} \textbf{239} (1996) \noopsort.

\bibitem{weibel:introduction}
Charles~A. Weibel: \emph{Introduction to Homological Algebra}, Cambridge
  studies in advanced mathematics 38, Cambridge University Press, Cambridge,
  1994.

\end{thebibliography}

\end{document}